\documentstyle[12pt]{article}
\def\ba{\begin{array}}
\def\ea{\end{array}}	
\def\be{\begin{equation}}
\def\ee{\end{equation}}
\def\bem{\begin{em}}
\def\eem{\end{em}}
\def\ot{\otimes}

\def\c{\times}
\def\a{\alpha}
\def\b{\beta}
\def\g{\gamma}

\def\lb{\langle}
\def\rb{\rangle}
\def\ep{\epsilon}
\def\ra{\longrightarrow}

\def\cb{{\cal B}}

\def\ce{{\cal E}}
\def\cf{{\cal F}}
\def\cg{{\cal G}}
\def\ch{{\cal H}}
\def\cm{{\cal M}}
\def\cn{{\cal N}}
\def\co{{\cal O}}

\def\cR{{\cal R}}

\def\cu{{\cal U}}
\def\cv{{\cal V}}
\def\cw{{\cal W}}
\def\cx{{\cal X}}

\def\tg{\triangle}
\def\lb{\langle}
\def\rb{\rangle}
\def\da{\downarrow}
\def\ua{\uparrow}
\def\ul#1,{\underline{#1}}
\newfont{\numb}{msbm10}
\def\com{\mbox{\numb C}}

\def\pom{\mbox{\numb P}}

\def\ps#1,#2,{\Psi_{#1,#2}}
\def\id#1,{id_{#1}}
\begin{document}
\title{SOME COMMENTS ON CATEGORIES, PARTICLES AND INTERACTIONS}
\author{Wladyslaw  Marcinek
\footnotemark[1]\\
Universit\"at Kaiserslautern, FB Physik\\
D-67663 Kaiserslautern, Germany} 
\date{}
\maketitle
\footnotetext[1]{
Permanent address: Institute of Theoretical Physics, University 
of Wroclaw, Plac Maksa Borna 9, 50-204 Wroclaw, Poland\\
e-mail: wmar@ift.uni.wroc.pl}
\date{}
\maketitle
\begin{abstract}
An algebraic formalism for the study of interacting particle systems 
is considered. Particle processes are described in terms of the
category theory. The problem for the unique description of these 
processes is discussed. Categories relevant for this subject are 
described. The concept of generalized transmutations of interacting
particle systems is introduced. The connection with a system with 
some generalized statistics is explained. 
\end{abstract} 
\newpage
\section{Introduction}
It is well--known that particles like baryons, nuclei, atoms or molecules 
are characterized by their own specific excitation spectrum. The existence 
of these spectra is one of the fundamental properties of the structure of 
matter. Suppose that we have a particle system with a collection of 
bound--state energy levels. There is a ground state and there are several 
excited states. It is also well--known that there are some transitions 
between different levels. They are results of interactions of the system 
with an external field or with other particle systems. These transitions 
which need more energy are known as excitation processes. They are the 
result of absorption of quanta of an external field.  On the other hand 
there are transitions connected with an energy spending, they correspond 
to some decay processes. It is interesting that there are also processes 
with no energy change. Let us  consider for instance these processes which 
can be described as sequences of vertex interactions of particle charges 
with an external quantum field. Charged particles are transformed under 
these interactions into a composite nonlocal discrete system which 
contains charges and quanta of the field. These systems are said to be 
{\it dressed particles} \cite{qsym}. We describe the structure of these 
dressed particles as a lattice with $n$ sites, $n = 1, 2, \ldots$. 
Every lattice site is a center for a vertex interaction of charge 
with the external quantum field. We assume that there are $N$ elementary
excitation states on every lattice sites. There are also collective
excitation states. It should be interesting to consider the general 
formalism for the study of all possible collective excitations of 
such system. We can imagine our lattice as a $d$-dimensional space 
(a manifold) equipped with $n$ distinct points as lattice sites.
One can consider a quantum dot, spin chains, or a set of vertex 
interactions of particles moving in two--dimensional space under 
influence of transversal magnetic field as examples of such systems 
corresponding for $d=0,1,2$, respectively. Our fundamental assumptions 
are that a collection of

$\bullet$ initial configurations of the system, 

$\bullet$ elementary particle processes.\\
is given. It seems to be natural to assume that every possible 
configurations of the system can be obtained as a result of certain
physical processes. We also assume that every process can be described 
as a sequence of elementary ones. These elementary processes 
represent elementary acts of lattice interactions. If all final 
configurations for the system under consideration can be described 
in an unique way as a result of transformation of an initial 
configuration, then we say that the system is equipped with a 
{\it category symmetry} or {\it coherent evolution}. 
Every such transformation is said to be a 
{\it evolution transformation}. This means that our 
category symmetry is in fact a formalism for the description
of particle interactions. The problem is to determine for a given 
system the smallest collection of symmetry transformations generating 
all others in an unique way.

The classical notion of the concept of symmetry in physics is based 
on group theory. The role played by the group representation 
theory for the study of symmetries in particle physics is well--known. 
The construction of a tensor product of representations is essential 
for such study. For instance, it allows to built states for composite
systems of particles from single particles ones. The unitary symmetry
and corresponding quark model is here a good example. It is known that 
we need a comultiplication for a tensor product of representations. It
is interesting that a comultiplication does exist for a large class
of $q$--deformed universal enveloping algebras. Hence they provide new 
possibilities for the study of particles, fields and their interactions
in mathematical physics. Also categories which contain a bifunctor 
$\ot :\cm\times\cm\ra\cm$ called a monoidal operation \cite{ML,jst,jet} 
provide some additional possibilities. Such categories are said to be 
monoidal. They can contain more structures and operations. The name
"monoidal" indicate, that there is one essential operation, just the 
monoidal  one, other operations play an auxiliary role. The bifunctor 
$\ot$ plays a similar role like the usual tensor product of group 
representations. 

Note that a related subject has been studied previously by several 
authors \cite{Lub,Gu2,man} and others. Formalism corresponding to 
the braided symmetry has been developed mainly by Majid 
\cite{SM,Maj,bm,sma,qm,smc}.
Categories in the context of quantum groups has been presented by
Kassel \cite{kas}. The application of categories in the topological 
quantum field theory has been considered by Sawin \cite{saw}. Similar 
formalism has been also developed previously by the author \cite{qsym}. 
Note that all these studies can be included in our general scheme. 
One can also consider the $q$--extended supersymmetry concept 
\cite{iluz,chung} as a particular example of our general formalism. 
Our considerations are mainly motivated by applications 
for the investigation of interacting systems of particles in 
low--dimensional spaces, but there are more different possible 
applications. It is known that the study of certain integrable 
models on a lattice leads to the investigation of some new formalism 
\cite{bus}. Hence it is interesting to study all these additional 
possibilities for the developing of the formalism beyond of the 
quantum mechanics and field theory. 

In this paper we are going to study these additional possibilities 
in a general manner in terms of monoidal categories. All our 
considerations are on abstract algebraic level. We would 
like to consider the most fundamental algebraic structures 
suitable for the description of particle interactions. 
We would like also to discuss the physical application for the  
classification of interacting particle systems. One can further 
develop our concept in terms of quantum von--Neumann algebras and 
their representations \cite{yam}. The paper is organized as follows. 
In Section 2 the general concept of particle interactions
is considered in terms of monoidal categories. 
Particle processes are described as certain transformations
of categories. The essential problem is the unique description, 
it is related to the coherence in categories \cite{ML}. 
In Section 3 categories relevant for our goal are described in details. 
Commutation relations for creation and annihilation quantum processes 
are described as certain  specific transformations. 
They lead to the system with generalized statistics \cite{gsi}. 
An introduction to the category theory is given in the Appendix. 
We believe that our approach can be 
useful for the deeper understanding of such new methods in 
quantum optics or both condensed matter and particle physics.
\section{General considerations}
Let us consider a system of hard core particles moving on certain
$d+1$--dimensional space--time manifold under influence of some
external field. All our considerations are based on the assumption 
that there is the vacuum state ${\bf 1}\equiv |0\rb$, the lowest 
energy elementary excited states $\{|i\rb\}_{i=1}^N$, and their
conjugated states ${\bf 1}^{\ast}\equiv\lb 0|$, and 
$\{\lb i|\}_{i=N}^1$ with the scalar product $\lb i|j\rb\in\com$.
We also assume that there are collective excitations of the system which 
can be described as a result of certain multiple product of elementary 
excitations. There is an energy gap between the vacuum state 
${\bf 1}\equiv |0\rb$ and the lowest energy excited state $|i\rb$. 
A finite set of $N$ operators $L = \{x^i\}_{i=1}^N$ is given 
as a starting point for our considerations. Every such operator  
act on the Hilbert space of functions on $d$--dimensional space.
We also assume that these operators transform functions representing 
the ground state of the system into states representing elementary 
excitations of the system, i.e. $|i\rb := x^i |0\rb$. In this way 
these operators represent elementary excited states of our system. 
Hence these operators are said to be primary. For the description of 
other states representing for instance collective excitations we need 
a product of operators. Such product need not forms a closed algebra 
but it must be defined in an unique way. A product of $n$ arbitrary 
primary operators should represents a collective $n$--tuple excitation.
It is an analog of a $n$ particle state. We would like to study 
such product in terms of the category theory. If $x^i$ is a primary 
operator, then there is a corresponding vector space $\cu = \cu (x^i)$. 
It is formally a $\com$--linear span of $x^i$, i. e. 
$\cu (x^i) := \{\a x^i;\a\in\com\}$. The $\com$--linear span of the 
ground state ${\bf 1}$ is denoted by $I$. It said to be the unit 
object. If $\cu$ and $\cv$ are  
$\com$--linear spans of $x^i$ and $x^j$, respectively, then the
linear span corresponding for certain product of these operators
is denoted by $\cu\ot\cv$ and is also said to be a product of $\cu$ 
and $\cv$. If for example $\cu$ represents charged particles 
excitation and $\cv$ some quanta, then 
the product $\cu\ot\cv$ describes the composite system containing 
both particles and quanta. This means that the operation 
$\ot : \cu\times\cv\ra\cu\ot\cv$ describes the "composition" process
of states. Such process tell us how to built a space of composite 
quantum states of the system from elementary ones. Hence it can be 
also understood as a generalization of the usual tensor product of 
group representations.  Observe that the arrow $\cu\ra\cu\ot\cv$ 
describes the process of absorption and the arrow $\cu\ot\cv\ra\cu$ 
describes the process of emission. 

Let us denote by $\pom$ the collection of all formal linear spans 
of primary operators, i.e. 
$\pom := \{\cu = \cu (x^i) : i=1,\cdots , N\}$. The collection of
complex conjugated spaces is denoted by $\pom^{\ast}$. 
We assume that an arbitrary sequence consists of the unit object $I$ 
or spaces from the collection $\pom$ or $\pom^{\ast}$ represents 
initial configuration of our system. These configurations can be
transformed into some new ones by a set of certain transformations.
These transformations represent certain physical processes like
composition, emission, absorption, etc... 
It is obvious that these transformations can be coherent or not. 
Coherence for a set of transformations means path--independent 
construction of these transformations.  Note that the coherence 
problem can be expressed graphically in terms of tangle tree 
operads \cite{moz}. Our goal is the construction 
of a collection of transformations which transform in an unique
way initial configurations into final ones -- representing the 
result of interactions. We denote $\ce$ the generating set for
these transformations. Let us consider some examples.\\
{\bf Example 1.} If $\ce$ contains only one operation, namely
the product $\ot :\pom\c\pom\ra\pom\ot\pom$, then starting from 
this product we can construct a set of multiproducts 
$\ot^n : \pom^{\times n}\ra\pom^{\ot^n}$ such that 
$\ot^2\equiv\ot$ and every multiproduct $\ot^n$ for $n>2$ 
can be calculated by an iteration procedure. 
Such procedure need not be unique. For instance for $n = 3$ 
we obtain $\ot^3 := \ot\circ(\ot\times id)$, But we also obtain
$\ot^3 := \ot\circ(id\times \ot)$. Hence for the uniqueness we need
some additional assumptions like the associativity
constraints, see the Appendix for more details.\\
{\bf Example 2.} For the (left) $\ast$--operation we use the standard
relations
\be
\cu^{\ast\ast} = \cu, \quad (\cu\ot\cv)^{\ast}=\cv^{\ast}\ot\cu^{\ast}.
\ee
In this case $\ce := \{\ot . \ast\}$.\\
{\bf Example 3.} We introduce a generating set 
$g(\pom) := \{g_{\cu} : \cu\in\pom\}$ of $I$--valued mappings 
$g_{\cu} : \cu^{\ast}\ot\cu\ra I$, where
\be
g_{\cu} (x^{\ast i}\ot x^j)\equiv(x^{\ast i}|x^j) := \lb i|j\rb
\ee
for pairing $g$. The extension $g_{\cu\ot\cv}$ of this pairing to 
the product $\cu\ot\cv$ is a problem. We need here the following
commutative diagram
\be
\ba{ccc}
&id_{\cv^*} \ot g_{\cu} \ot id_{\cv}&\\
\cv^* \ot \cu^* \ot \cu \ot \cv& \ra&\cv^* \ot \cv\\
&&\\
\parallel&&\downarrow g_{\cv}\\
&&\\
(\cu \ot \cv)^* \ot (\cu \ot \cv)&\ra&I\\
&g_{\cu \ot \cv}&
\ea
\ee
for the extension. We can introduce the set
$g(\pom^{\ast}) := \{g_{\cu} : \cu\in\pom\}$ of $I$--valued 
mappings $g_{\cu^{\ast}} : \cu\ot\cu^{\ast}\ra I$ in a similar way.
For the extension we use the diagram
\be
\ba{ccc}
&id_{\cu} \ot g_{\cv^{\ast}}\ot id_{\cu^{\ast}}&\\
\cu\ot\cv\ot\cv^{\ast}\ot\cu^{\ast}&\ra&\cu\ot\cu^{\ast}\\
&&\\
\parallel&&\downarrow g_{\cu^{\ast}}\\
&&\\
(\cu\ot\cv)\ot(\cu\ot\cv)^{\ast}&\ra&I\\
&g_{(\cu\ot\cv)^{\ast}}&
\ea
\ee
In the way we can construct a category $\cm := \cm (\pom , \ce)$ 
whose objects are multiple products of object of $\pom$ and morphisms 
are obtained by iteration procedures applied to operations from the 
set $\ce$. The monoidal operation of $\cm$ can also be obtained by 
the proper iteration of the initial product $\ot$. 

Note that there is the uniqueness problem with the construction 
of the category $\cm = \cm (\pom , \ce)$ over $\pom$. One can 
construct many different categories $\cm$ for a given initial
collection of spaces $\pom$ with different generating set $\ce$. 
We denote by ${\bf Cat}(\pom)$ the class of all these categories.
Let $\cm (\pom , \ce)$ and $\cn (\pom, \ce')$ be two such categories.
An arbitrary functor $\cf : \cm\ra\cn$ which transform the set of
operations $\ce$ into $\ce'$ is said to be a {\it generalized
transmutation}. In this case we say that the set $\ce$ is
transmuted into $\ce'$. The category $\cn$ is then said to
be {\it functored} over $\cm$ \cite{smaj}.

If $\ce := \{I , \ot\}$ and $\ce' := \{I', \ul \ot,\}$, then the 
corresponding generalized transmutation $\cf : \cm\ra\cn$ is a 
monoidal functor of categories $\cm$ and $\cn$. This means that
it is a triple
$$
\cf := \{\cf,\varphi_2,\varphi_0\} : \cm \ra \cn
$$
which consists of a functor $\cf : \cm\ra\cn$,
a natural isomorphism
$$
\varphi := \varphi_{2,\cu,\cv} : 
\cf\cu\ul \ot, \cf\cv \ra \cf(\cu \ot \cv),
$$
and an isomorphism
$
\varphi_0 : I \ra \cf I = I'
$, 
such that the following diagrams
\be
\ba{ccc}
&\varphi_2\ul \ot, id& \\
\cf\cu\ul \ot, \cf\cv\ul \ot, \cf\cv& \longrightarrow&
\cf(\cu \ot \cv)\ul \ot, \cf\cw\\
&&\\
id\ot\varphi_2\downarrow&&\downarrow\varphi_2\\
&&\\
\cf\cu\ul \ot, \cf(\cv \ot \cw)&\longrightarrow&
\cf(\cu\ot\cv\ot\cw)\\
&\varphi_2&
\ea
\ee

\vspace{0.5cm}

\be
\ba{rcc}
&\varphi_2& \\
\cf I\ul \ot, \cf\cu&\longrightarrow&\cf (I\ot\cu)\\
&&\\
\varphi_0\ul \ot, id\uparrow&\swarrow&\cf (l_{\cu})\\
&&\\
\cf\cu&&
\ea
\ee

\vspace{0.5cm}

\be
\ba{rcl}
&\varphi_2& \\
\cf\cu\ul \ot, \cf I&\longrightarrow&\cf(\cu\ot I)\\
&&\\
id\ul \ot, \varphi_0\uparrow&\swarrow&\cf (r_{\cu})\\
&&\\
\cf\cu&&
\ea
\ee
are commutative. If $\varphi_2$ and $\varphi_0$ are identities, then 
$\cf$ is said to be strict. For the $\ast$--operation and pairing we 
have
\be
(\cf (\cu))^{\ast} = \cf (\cu^{\ast}), 
\quad g_{\cu} = g'_{\cf (\cu)};
\ee
respectively. A generalized transmutation $\cf : \cm\ra\cn$ 
is said to be cross symmetric if the following diagram
\be
\ba{ccc}
&\varphi_2&\\
\cf\cu^{\ast}\ul \ot, \cf\cv&\longrightarrow&\cf(\cu^{\ast}\ot\cv)\\
&&\\
\Psi'\downarrow&&\downarrow \cf(\Psi)\\
&&\\
\cf\cv\ul \ot, \cf\cu^{\ast}&\longrightarrow&
\cf (\cv\ot\cu^{\ast})\\
&\varphi_2&\\
\ea
\ee
is commutative for every generating objects $\cu, \cv$ of $\cm$, 
where $\Psi$ and $\Psi'$ are the cross symmetries in $\cm$ and $\cn$,
respectively. 
In the case of a braid symmetries we have the following diagram
\be
\ba{ccc}
&\varphi_2&\\
\cf\cu\ul \ot, \cf\cv&\longrightarrow&\cf(\cu\ot\cv)\\
&&\\
\Psi'\downarrow&&\downarrow \cf(\Psi)\\
&&\\
\cf\cv\ul \ot, \cf\cu&\longrightarrow&\cf(\cv\ot\cu)\\
&\varphi_2&\\
\ea
\ee
Note that the functor $\cf : \cm\ra\cn$ generalizes the Majid
concept of transmutation of braid statistics \cite{mat}. If
$\cm$ is the category of cobordisms of smooth manifolds and
$\cn$ is a category of vector spaces, that the generalized
transmutation is known as the Topological Quantum Field Theory
\cite{saw}.\\
{\bf Example 4.} Let $H$ and $H'$ be two Hopf algebras. If 
$h : H\ra H'$ is a Hopf algebra homomorphism, then we can
introduce the transmutation $\cf : \cm^{H}\ra\cm^{H'}$ of
categories of right comodules as follows. The functor $\cf$
acts as the identity functor on arbitrary comodule $\cu$ but
the coaction $\rho_{\cu}:\cu\ra\cu\ot H$ transform into new
one, namely into $\rho'_{\cu}:\cu\ra\cu\ot H'$, where
\be
\rho'_{\cu} := (\id \cu, \ot h)\circ\rho_{\cu}.
\ee
For coquasitriangular Hopf algebras $H$ and $H'$ with 
coquasitriangular structures $\lb -,-\rb :H\ot H\ra I$ and 
$\lb -,-\rb' :H'\ot H'\ra I$, respectively, we obtain 
\be
\ba{c}
\Psi'_{\cu,\cv}(u\ot v) = \Sigma\lb h(v_1),h(u_1)\rb' \ v_0 \ot u_0 ,
\label{cons}
\ea
\ee
where $\rho (u) = \Sigma u_0\ot u_1\in\cu\ot H$, 
$\rho (v) = \Sigma v_0\ot v_1\in\cv\ot H$ for every 
$u\in\cu , v\in\cv$, and  
$\lb k, l\rb = \lb h(k), h(l)\rb'$ for every $k,l\in H$.\\
{\bf Example 5.} Let $H := \com G$ and $H' := \com G'$ be group
algebras, where $G$ and $G'$ are Abelian groups equipped with
factors $\ep$ and $\ep'$, respectively. Then the transmutation
$\sf : \cm (G, \ep)\ra\cm (G', \ep')$ is determined by a group 
homomorphism $h : G\ra G'$ such that
\be
\ep (\a , \b) = \ep (h(\a) , h(\b))
\ee
for $\a , \b\in G$.
\section{Commutation relations}
Let us denote by $\pom (n)$ the collection of all tensor products 
of the form
\be
\ba{l}
\cu_{i_1}\ot\cdots\ot\cu_{i_n},
\ea
\ee
for all $\cu_{i_1},\cdots, \cu_{i_n}\in\pom$. 
We introduce $\pom^{\ast}(n)$ in a similar way. We also
introduce the collection $(\pom^{\ast}\ot\pom)(n)$
of sequences of the form
\be
\ba{l}
\cu^{\ast}_{j_n}\ot\cdots\ot\cu^{\ast}_{j_1}\ot
\cu_{i_1}\ot\cdots\ot\cu_{i_n}.
\ea
\ee
The collection $(\pom\ot\pom^{\ast})(n)$ can be defined in an 
obvious way. We have the following examples.\\
{\bf Example 4.} Let us denote by 
$\cm := \cm (\pom , I,  \ot , \ast , g)$ the category 
$\cm := \cm (\pom , \ce)$, where $\ce := \{I, \ot , \ast , g\}$.
We introduce two sets of transformations 
$a^+ := \{a^+_{\cu}:\pom (n)\ra\pom (n+1)\}$ 
and $a^- := \{a^-_{\cu^{\ast}}:\pom (n)\ra\pom (n-1)\}$, where 
\be
\ba{c}
a^+_{\cu}(\cu_{i_1}\ot\cdots\ot\cu_{i_n}):= 
\cu\ot\cu_{i_1}\ot\cdots\ot\cu_{i_n} ,
\ea
\ee
and
\be
\ba{c}
a^-_{\cu^{\ast}}(
\cu_{i_1}\ot\cdots\ot\cu_{i_n}):= 
g_{\cu}(\cu^{\ast}_{j_1}\ot\cu_{i_1})\ot\cdots\ot\cu_{i_n}
\ea
\ee
where
\be
g_{\cu}(\cu^{\ast}_{j_1}\ot\cu_{i_1})
\left\{
\ba{cl}
g_{\cu}&\mbox{for}\quad\cu^{\ast}_{j_i}\equiv\cu^{\ast},
\cu_{i_1}\equiv\cu\\
0&\mbox{otherwise}
\ea
\right. ,
\ee
for $\cu_{i_1},\cdots, \cu_{i_n},\cu$ and 
$\cu^{\ast}_{j_1},\cu^{\ast}\in\pom$. Here $\cu^{\ast}$ 
represents a quasihole and $\cu$ -- a quasiparticle.
Two different objects $\cu ,\cv\in\pom$ represent (quasi-)
particle states of two different sorts. There are no identical
particles. It is easy to see that we have the following set of 
relations
\be
a^-_{\cu^{\ast}}\circ a^+_{\cu} = g_{\cu}{\bf 1},
\ee
where $\cu\in\pom$. These relations are in fact the commutation 
relations for the system equipped with the infinite statistics 
\cite{gsd}. One can use the relation
$$
\hat{a}_{\cu}\circ \hat{a}_{\cu} = \hat{a}_{\cu\ot\cv},
$$
where $\hat{a}_{\cu}$ stands for $a^+_{\cu}$ or $a_{\cu^{\ast}}$,
for the extension of commutation relations corresponding for
monoidal products of generating objects. These relations seems to be
simple, but they lead to well--defined operator algebras \cite{braj}.  
Observe that we have here elementary quantum processes of two sorts,
namely creation and annihilation.\\
{\bf Example 5.} We assume that $\ce$ contains the cross
symmetry $\Psi_{cross}$ in addition to the previous example.
The corresponding category is denoted by
$\cm_{cross} := \cm (\pom , I,  \ot , \ast , g , \Psi_{cross})$.
We need here a collection
$\Psi (\pom) := \{\Psi_{\cu^{\ast},\cv} : \cu^{\ast}, \cv\in\pom\}$ 
as the initial data for the description of exchange processes of 
(quasi-) particles and (quasi-) holes, see the Appendix. 
We have here the following relations 
\be
\ba{c}
b^-_{\cu^{\ast}}\circ b^+_{\cu} - \ b^+_{\cu}\circ b^-_{\cu^{\ast}} 
\circ \ps \cu^{\ast}, \cu, := g_{\cu} {\bf 1},
\label{cca}
\ea
\ee
where 
\be
\ba{l}
b^+_{\cu} := a^+_{\cu},\\
b^-_{\cu^{\ast}}(\cv_1\ot\cv_2 ) := \\ 
\left[
(a^-\ot id_{\cv_2}) - (id_{\cv_1 }\ot a^-)\circ 
(\ps \cu^{\ast}, \cv, \ot id_{\cv_2 }) 
\right] 
(U^{\ast}\ot\cv_1\ot\cv_2 ),
\ea
\ee
and
\be
a^- (\cu^{\ast}\ot\cv) := a^-_{\cu^{\ast}}(\cv).
\ee
In this way we have here a collection of elementary processes $\ce$ 
which contains creation, annihilation and exchange processes.\\
{\bf Example 6.} We replace the cross symmetry by the braid one.
In this case we obtain the following relations
\be
\ba{c}
c^-_{\cu^{\ast}}\circ c^+_{\cu} - \ c^+_{\cu}\circ c^-_{\cu^{\ast}} 
\circ \ps \cu^{\ast}, \cu, := g_{\cu} {\bf 1},
\label{bcr}
\ea
\ee
and in addition 
\be
\ba{l}
c^-_{\cu^{\ast}}\circ c^+_{\cv} - \ c^+_{\cv}\circ c^-_{\cu^{\ast}} 
\circ \ps \cu^{\ast}, \cv, = 0,\\
c^+_{\cu}\circ c^+_{\cv} - \ c^+_{\cv}\circ c^+_{\cu} 
\circ \ps \cu, \cv, = 0 ,\\
c^-_{\cu^{\ast}}\circ c^-_{\cv^{\ast}} - \ 
c^-_{\cv^{\ast}}\circ c^-_{\cu^{\ast}} 
\circ \ps \cu^{\ast}, \cv^{\ast}, = 0 .
\label{ccd}
\ea
\ee
Note that for the braid symmetry there are additional elementary 
quantum processes, namely the exchange processes of identical 
particles on lattice in two dimensional case.\\ 
{\bf Example 7.} Let
$\cm_{cross} := \cm (\pom , I,  \ot , \ast , g , \Psi_{cross})$
be a category with commutation relation1s like in the Example 5.
we denote by 
$\cn_{cross}:=\cn (\pom,I, \ul \ot,,\-\ast , g',\Psi'_{cross})$
a second category with a new cross $\Psi'_{cross}$ and pairing 
$g'$. One can define the following two sets 
$c^+ := \{c^+_{\cu}:\pom (n)\ra\pom (n+1)\}$ 
and $c^- := \{c^-_{\cu^{\ast}}:\pom (n)\ra\pom (n-1)\}$
of operators in it. For a cross symmetric generalized transmutation 
$\cf : \cm\ra\cn$, we have the relation for these operators
\be
\ba{c}
c^-_{\cf(\cu^{\ast})}\circ c^+_{\cf(\cu)} -
 \ c^+_{\cf(\cu)}\circ c^-_{\cf (\cu^{\ast})} 
\circ \Psi'_{\cf(\cu^{\ast}),\cf (\cu)} := g'_{\cf(\cu)} {\bf 1},
\label{ccc}
\ea
\ee
Note that the category $\cn$ can be braided or symmetric.
In these cases we obtain additional relations such as (\ref{ccc}).

It is obvious that the concept of category symmetries is related 
to the systems with generalized statistics \cite{gsi,gsd}. Note 
that the braid commutation relations, consistency conditions and 
corresponding Fock space representation with well--defined 
scalar product has been considered previously, see \cite{quon} for 
instance. Some interesting examples of related formalism 
has been studied previously by Fiore \cite{fio}. 
Observe that the above concept of category symmetries can be 
futher developed in a few respects. One can consider the corresponding 
noncommmutative calculi, It should be interesting to study 
Hamiltonians in terms of described here creation and annihilation 
operators and study the concrete physical models. 
\section*{Appendix}
Let us briefly recall the fundamental concept of the category 
theory for the fixing of notation. For more details see the
textbook of Mac Lane \cite{ML}.
A category $\cm$ contains a collection $\co b(\cm)$ of objects
and a collection $hom(\cm)$ of arrows (morphisms).
The collection $hom(\cm)$ is the union of mutually disjoint 
sets $hom(\cu, \cv)$ of arrows $f : \cu \ra \cv$ from $\cu$ to $\cv$ 
defined for every pair of objects $\cu, \cv \in \co b(\cm)$. It may 
happen that for a pair $\cu, \cv\in\co b(\cm)$ the set $hom(\cu, \cv)$ 
is empty. The associative composition of morphisms is also defined.
A functor $\cf : \cm \ra \cn$ of the category $\cm$ into the category 
$\cn$ is a map which sends objects of $\cm$ into objects of $\cn$ and 
morphisms of $\cm$ into morphisms of $\cn$ such that
$\cf(f \circ g) = \cf(f) \circ \cf(g)$
for every morphisms $f : \cv \ra \cw$ and $g : \cu \ra \cv$
of $\cm$. The generalization to multifunctors is obvious.
One can consider an arbitrary object of a category as an example 
of constant functor. For instance an $n$--ary functor 
$\cf : \cm^{\times n}\ra\cn$
sends an $n$--tuple of objects of $\cm$ into an object of $\cn$.
The corresponding condition for morphisms is evident.
In this paper we restrict our attention for a description how
functors act on objects, we omit the action on morphisms 
for simplicity. The reader can complete our description.

Now we recall the concept of natural transformations.
Let $\cf$ and $\cg$ be two functors of the category $\cm$ into
the category $\cn$. A natural transformation $s : \cf \ra \cg$
of $\cf$ into $\cg$ is a collection of morphisms
$s = \{s_{\cu} : \cf (\cu) \ra \cg (\cu), \cu \in {\cal O}b(\cm)\}$
such that
\be
s_{\cv} \circ \cf (f) = \cg (f) \circ s_{\cu}
\ee
for every morphism $f : \cu \ra \cv $ of $\cm$. The set of all
natural transformations of $\cf$ into $\cg$ is denoted by
$\cn at(\cf,\cg)$. It is easy to see that the composition
$t \circ f$ of natural transformation $s$ of $\cf$ into $\cg$
and $t$ of $\cg$ into $\ch$ is a natural transformation of
$\cf$ into $\ch$. If $\cf \equiv \cg$, then we say that the natural
transformation $s : \cf \ra \cg$ is a natural transformation
of $\cf$ into itself.

Now let us briefly explain the notions of monoidal categories 
\cite{ML,jst} adopted for our goal. 
{\it A monoidal category} $\cm\equiv\cm(\ot, I)$ is in fact 
a category $\cm$ equipped with a monoidal operation (a bifunctor) 
$\ot : \cm\times\cm\ra\cm$, a unit object $I$, and 
collections of natural isomorphisms:\\
(i) an associativity constraint $\psi = \{\psi_{\cu,\cv,\cw}:
(\cu\ot\cv)\ot\cw\ra\cu\ot (\cv\ot\cw)$\},\\
(ii) a left unity constraint $l = \{l_{\cu}:I\ot\cu\ra\cu\}$\\
(iii) and a right unity constraint $r = \{r_{\cu}:\cu\ot k\ra\cu\}$\\
such that the following diagrams
\be
\ba{ccc}
&(\cu \ot \cv)\ot(\cw \ot \cx)&\\
\psi_{\cu\ot\cv,\cw,\cx}\nearrow&&\searrow\psi_{\cu,\cv,\cw\ot\cx}\\
&&\\
((\cu\ot\cv)\ot\cw)\ot\cx&&\cu\ot(\cv\ot(\cw\ot\cx))\\
&&\\
\psi_{\cu,\cv,\cw}\ot id\da&&\ua id \ot \psi_{\cv,\cw,\cx}\\
&&\\
(\cu\ot(\cv\ot\cw))\ot\cx&\ra&\cu\ot((\cv\ot\cw)\ot\cx)\\
&\psi_{\cu,\cv \ot \cw,\cx}&
\ea
\ee
\be
\ba{ccc}
&\psi_{\cv ,I,\cw}&\\
(\cv \ot I) \ot \cw&\ra&\cv \ot (I \ot \cw)\\
&&\\
r_{\cv} \ot id \searrow&&\swarrow id \ot l_{\cw}\\
&&\\
&\cv \ot \cw&
\ea
\ee
commute. It is interesting that in a monoidal category any diagram built 
from the constraints $\psi , l, r$, and the identities by composing and 
tensoring, commutes. This is just the famous Mac Lane's coherence
theorem. A monoidal category $\cm$ is said to be {\it strict}, if all 
natural isomorphisms $\psi_{\cu,\cv,\cw},l_{\cu},r_{\cu}$ are identity. 
It is also interesting that every monoidal category is equivalent to
certain strict one. This means that we can restrict our attention to
strict monoidal categories.

{\it A (left) $\ast $-operation} in a monoidal category $\cm$ 
is a transformation $(-)^{\ast}$ of functor $\ot$ into the opposite
functor $\ot^{op}$ such that
\be
\ba{cc}
(-)^{\ast \ast} = id_{\cm},
&(-)^{\ast}\circ\ot = \ot^{op}\circ (-)^{\ast} 
\ea
\ee
where $\cu$ and $\cv$ are arbitrary objects of the category $\cm$.
{\it A (left) pairing} $g$ in the category $\cm$ is a transformation 
of the functor $(-)^{\ast}\ot -$ into $I$, where $I$ is a field
satisfying some compatibility axioms, see \cite{castat,bes}.
This means that $g$ is a set $g\equiv\{g_{\cu}\}$ of $I$--valued 
mappings 
\be
\ba{c}
g\equiv\{g_{\cu} : \cu^{\ast}\ot \cu \ra I, \cu \in\co b(\cm)\}
\label{parc}
\ea
\ee
Let $\cm$ be a monoidal category equipped with a (left) 
$\ast$-operation $(-)^{\ast}$ and a (left) pairing $g$, 
then such category is said to be a {\it category with (left) duality}
and it is denoted by 
$\cm = \cm_{left}(\ot ,\oplus ,I, \ast , g)$.
One can introduce a (right) duality structure in the 
category $\cm$ in a similar way. Note that both dualities in $\cm$
the right and the left one are in general two independent structures. 
But it is possible to introduce an additional structure which making 
these two structures equivalent. Such equivalence can be established 
by the following set of natural isomorphisms 
\be
\Psi \equiv \{\ps \cu^{\ast},\cv, : 
\cu^{\ast} \ot \cv \ra \cv \ot \cu^{\ast} \}.
\ee
where
\be
\ba{l}
\ps \cu^{\ast} \ot \cv^{\ast},\cw, 
= (\ps \cu^{\ast},\cw, \ot \id \cv,) 
\circ (\id \cu, \ot \ps \cv^{\ast},\cw,),\\
\ps \cu^{\ast},\cv \ot \cw,
= (\id \cv, \ot \ps \cu^{\ast},\cw,) 
\circ (\ps \cu^{\ast},\cv, \ot \id \cw,),
\label{heb}
\ea
\ee
for every objects $\cu ,\cv ,\cw$ in $\cm$. These transformations
are called a {\it generalized cross symmetry}, \cite{castat}.
We can identify the right and left duality in the category equipped 
with such generalized cross symmetry. The monoidal category equipped
with such symmetry is denoted by 
$\cm = \cm (\ot ,\oplus ,I , \ast, g, \Psi_{cross})$. 

Note that the generalized cross symmetry is not a braid symmetry
in general. For the braid symmetry in the category with duality
we need additional transformations like 
\be
\Psi \equiv \{\ps \cu,\cv, : \cu \ot \cv \ra \cv \ot \cu\}
\ee
for arbitrary objects $\cu ,\cv$ in $\cm$ and
\be
\Psi \equiv 
\{
\ps \cu^{\ast},\cv^{\ast}, : 
\cu^{\ast} \ot \cv^{\ast} \ra \cv^{\ast} \ot \cu^{\ast} 
\}
\ee
for objects $\cu^{\ast} ,\cv^{\ast}$ in $\cm$. 
We need also some new commutative diagrams for all these 
transformations and pairings. In fact a family of natural 
isomorphisms
\be
\ba{c}
\Psi \equiv \{\Psi_{\cu \ot \cw} : \cu \ot \cw \ra \cw \ot \cu \}
\label{bs1}
\ea
\ee
such that we have the following relations
\be
\ba{l}
\ps \cu \ot \cv,\cw, = (\ps \cu,\cw, \ot \id \cv,) 
\circ (\id \cu, \ot \ps \cv,\cw,),\\
\ps \cu,\cv \ot \cw, = (\id \cv, \ot \ps \cu,\cw,)  
\circ (\ps \cu,\cv, \ot \id \cw, ),
\label{hex}
\ea
\ee
is said to be {\it a braiding or a braid symmetry} on $\cm$. The 
monoidal category with unique duality and braid symmetry is said 
to be {\it rigid} \cite{Maj,bm,sma,qm}. This category is denoted 
by $\cm = \cm (\ot ,\oplus ,I , \ast , g. \Psi_{braid})$. 
If in addition we have the relation
\be
\ba{c}
\Psi_{\cu,\cv}^2 = id_{\cu \ot \cv},
\label{kwad}
\ea
\ee
for every objects $\cu, \cv \in \cm$, then the set $S := \{\Psi_{\cu,\cv}\}$ 
is said to be a (vector) symmetry or tensor symmetry and the corresponding 
category $\cm$ is called {\it a symmetric monoidal or tensor category}, 
see \cite{Lub,man}. 

Let us consider some examples of monoidal categories which
can be useful for the study of category symmetry. 
The most simple example of a monoidal category is provided by 
the category $\cv ect(k)$ of vector spaces over a field $k$.
The monoidal operation in this category is defined by the usual tensor 
product of vector spaces. Another example is given by the 
category $\cv ect_G (k)$ of $G$--graded vector spaces, where $G$
is a grading group. In the supersymmetry the grading group is the
group of integer $Z_2$. For anyons we have $G\equiv Z_n$, where 
$n>2$, \cite{smc}.There is a category $_{\cb}\cm$ of all left 
$\cb$-modules, where $\cb$ is an unital and associative algebra. 
Observe that the usual tensor product $\cu\ot\cv$ of two left 
$\cb$--modules $\cu$ and $\cv$ is not a left $\cb$-module but a left
$\cb\ot\cb$-module! Hence this category is not a monoidal category. But 
it is easy to see that in the particular case when $\cb$ is a bialgebra, 
i.e. we have a comultiplication $\tg : \cb\ra\cb\ot\cb$ in $\cb$, the 
category $_{\cb}\cm$ is monoidal.  
For instance there is the category $\cR_G$ of finite dimensional
representations of compact matrix quantum group $G$, \cite{wor}.
There is also a category of Hopf modules or crossed modules 
\cite{bes}. Observe that there is also a category $\cm^{\ch}$ 
of right $\ch$--comodules, where $\ch$ is a Hopf algebra. The monoidal 
operation in $\cm^{\ch}$ is given as the following tensor product of 
$\ch$--comodules
\be
\rho_{\cu\ot\cv} = (id \ot m_{\ch}) \circ (id \ot \tau \ot id)
\circ (\rho_{\cu} \ot \rho_{\cv}),
\ee
where $\tau : \cu\ot\ch\ra\ch\ot\cu$ is the twist, 
$m_{\ch} : \ch\ot\ch\ra\ch$ is the multiplication in $H$.

As an example for a category with duality we can give a category 
$_{\cb}\cm$ of left $\cb$--modules. In this case the monoidal 
operation corresponds to the tensor product of representations, 
the $\ast$-operation corresponds to the 
contragradient representation and the generalized 
cross symmetry corresponds to the intertwiner between an arbitrary 
representation and its contragradient. Hence it is also called 
{\it a statistics operator}. Note that if there is a bialgebra $\cb$ 
such that the category $\cm$ is equivalent to the category of left 
modules (representations) over $\cb$, then the bialgebra is said to 
be {\it a generalized symmetry} or {\it (bi-)algebra symmetry} for 
the given physical system. One can describe states in the quantum 
Hall effect as a result of symmetry in such generalized sense 
\cite{oma}. The symmetry algebra for Klein--Gordon equation on 
quantum Minkowski space is considered in \cite{kli}.

Note that the category of 
representations of the so--called weak Hopf algebra is rigid \cite{bos}. 
Also the category of quantum compact matrix groups of Woronowicz is 
rigid \cite{wor,yam}. For a coquasitriangular Hopf algebra $H$ with 
a coquasitriangular structure $\lb -,-\rb : H\ot H\ra I$ we obtain
the category $\cm^H$ of right $H$-comodules which is also braided 
monoidal \cite{qsym}. The braid symmetry 
$\Psi\equiv\{\ps \cu, \cv, : \cu\ot\cv\ra\cv\ot\cu; \cu, \cv\in Ob\cm\}$
in $\cm$ is defined by the equation
\be
\ba{c}
\ps \cu, \cv, (u\ot v) = \Sigma\lb v_1 , u_1 \rb \ v_0 \ot u_0 ,
\label{coin}
\ea
\ee
where $\rho (u) = \Sigma u_0 \ot u_1 \in \cu \ot H$, and
$\rho (v) = \Sigma v_0 \ot v_1 \in \cv \ot H$ for every 
$u\in \cu , v \in\cv$.

Let $G$ be an arbitrary group, then the group algebra $H := \com G$
is a Hopf algebra for which the comultiplication, the counit,
and the antypode are given by the formulae
$$
\ba{cccc}
\tg (g) := g \ot g,&\eta(g) := 1,&S(g) 
:= g^{-1}&\mbox{for} \ g \in G.
\ea
$$
respectively. If $H\equiv\com G$, where $G$ is an Abelian 
group, then the coquasitriangular structure on $H$ is given
as a bicharacter on $G$ \cite{mon}. Note that for
Abelian groups we use the additive notation. A mapping 
$\ep : G\c G \ra\com\setminus \{0\}$ is said to 
be a {\it bicharacter} on $G$ if and only if we have 
the following relations
\be
\ep (\a, \b + \g) = \ep (\a, \b) \ep (\a , \g), \quad
\ep (\a + \b , \g) = \ep (\a , \g) \ep (\b , \g)
\ee
for $\a, \b, \g \in G$. If in addition
\be
\ep (\a, \b) \ep (\b, \a) = 1, 
\ee
for $\a, \b \in G$, then $\ep$ is said to be 
a {\it normalized bicharacter} or a commutation 
factor on $G$ \cite{sch}. The category $\cm^H$ of right
comodules, where $H := \com G$ for certain Abelian group $G$ 
and $\lb -,-\rb\equiv\ep (-,-)$ is a bicharacter like above 
is denoted by $\cm (G, \ep)$. Note that if $\cu$ is a $H$-comodule, 
where $H = \com G$, then $\cu$ is a $G$-graded vector space, i.e 
$\cu = \bigoplus\limits_{\a \in G} \ \cu_{\a}$. This means that a 
coaction of $H := \com G$ on $\cu$ is equivalent to $G$-gradation
of $\cu$.

\newpage
\noindent {\bf Acknowledgments}\\
The author would like to thank to Prof, W. R\"uhl for kind
invitation to Universit\"at Kaiserslautern and discussion,
to M. Greulach for any other help.  
The work is partially sponsored by DAAD, Bonn, Germany, and by
Polish Committee for Scientific Research (KBN) under Grant 2P03B130.12.

\end{document}